\documentclass[12pt]{article}
\usepackage{amsmath, amssymb, graphicx, geometry, setspace, cite}
\usepackage[utf8]{inputenc}
\usepackage{hyperref}
\usepackage{graphicx}
\usepackage{subcaption}
\newtheorem{theorem}{Theorem}

\title{Detecting transitions from steady states to chaos with gamma distribution}
\author{Haiyan Wang\\
School of Mathematical and Natural Sciences\\
Arizona State University\\
Phoenix, AZ 85069\\
haiyan.wang@asu.edu\\\\
Ying Wang\\
Department of Mathematics\\
University of Oklahoma\\ 
Norman, OK 73019-3103\\
wang@ou.edu
}
\date{}

\begin{document}
\maketitle

\section*{Abstract}
In this paper, we introduce a novel method to identify transitions from steady states to chaos in stochastic models, specifically focusing on the logistic and Ricker equations by leveraging the gamma distribution to describe the underlying population. We begin by showing that when the variance is sufficiently small, the stochastic equations converge to their deterministic counterparts. Our analysis reveals that the stochastic equations exhibit two distinct branches of the intrinsic growth rate, corresponding to alternative stable states characterized by higher and lower growth rates.

Notably, while the logistic model does not show a transition from a steady state to chaos, the Ricker model undergoes such a transition when the shape parameter of the gamma distribution is small. These findings not only enhance our understanding of the dynamic behavior in biological populations but also provide a robust framework for detecting chaos in complex systems.

\textbf{keywords:} Chaos, steady state, logistic difference equation, Ricker difference equation, gamma distribution

%
%

\section{Introduction}
Identifying chaotic states through nonlinear dynamics analysis is essential for practical applications, enabling better prediction and management of complex, unpredictable behaviors in various systems \cite{hilborn2000chaos,strogatz2018nonlinear,earlywarning2023}. It is well-known that the logistic map showed that as a control parameter \(r\) increases, a system transitions from steady states to periodic behavior and eventually to chaotic dynamics \cite{may1976simple,feigenbaum1980quantitative}. While bifurcation analysis played a foundational role in the discovery of chaos, the field has since expanded to encompass other methods, such as Lyapunov exponents \cite{wolf1985lyapunov} and recurrence quantification analysis \cite{eckmann1987recurrence}.

Chaos often emerges in difference equations, which model discrete-time systems. The simplest and most famous examples are the logistic map:
\begin{equation} \label{eq:dlogistic3}
x_{t+1} = r x_t ( 1 - x_t ),
\end{equation}
and the Ricker map:
\begin{equation}
x_{t+1} = x_t e^{r(1 - x_t)}
\label{dricker1}
\end{equation}
where $x_t$ represents the state of the system at time step $t$, and $r$ is the growth rate.  For \eqref{eq:dlogistic3}, chaos emerges at higher growth rates (\(r > 3\)) following a period-doubling route to chaos \cite{may1976simple, feigenbaum1980quantitative,li1975}. For \eqref{dricker1}, incorporates exponential density dependence, leading to sharper suppression of population growth at higher densities \cite{ricker1954stock}. This sharper feedback allows the Ricker model to exhibit chaotic dynamics at much lower growth rates (\(r \sim 2\)) compared to the logistic model  \cite{may1976simple}. 

The impact of noise depends on the system's structure, the type of noise introduced, and the underlying deterministic dynamics \cite{Allen,crutchfield1982fluctuations,kot2001elements,moss1989noise}. For example, in the stochastic logistic map,
\begin{equation}
    x_{t+1} = r x_t (1 - x_t) + \xi_t,
\end{equation}
where \(\xi_t\) is a random noise, often taken from a Gaussian distribution, plays a crucial role, such as population dynamics, chemical reactions, and financial markets. \cite{sato2018dynamical} provide a detailed characterization of these stochastic bifurcations by analyzing the Lyapunov exponent and the dichotomy spectrum. \cite{erguler2008statistical} studied how noise alters the probability density function (PDF) of trajectories in the stochastic logistic map. The stochastic logistic map has also been studied from the perspective of algorithmic probability and simplicity bias in \cite{hamzi2024simplicity}. \cite{li2008effects} showed that noise can increase the stability of certain periodic orbits by altering their probability densities.  In addition, \cite{Yan2024} introduces a modified Ricker model and evaluates the effects of variances and expectations of the discrete uniform random variable on the dynamics. These works suggest that the distribution of $x_t$ plays a significant role in the study of stochastic population dynamics. 

In this paper, we study stochastic logistic equation \eqref{eq:logistic3} and stochastic Ricker equation \eqref{ricker1} 
\begin{equation} \label{eq:logistic3}
X_{t+1} = r X_t \left( 1 - X_t \right) \epsilon_t,
\end{equation}
\begin{equation}
X_{t+1} = X_t e^{r(1 - X_t)} \epsilon_t
\label{ricker1}
\end{equation}
where \( X_t \) is a distribution of population size at time \( t \), \( r \) is the intrinsic growth rate,  $\epsilon_t$ is a small nonnegative perturbation distribution representing stochastic effects, assumed it is independent of $X_t$ and its mean is $1$ ($E[\epsilon_t] = 1$). We investigate the relation between the stochastic equations \eqref{eq:logistic3} and \eqref{ricker1}, and their deterministic versions \eqref{eq:logistic3} and \eqref{ricker1} respectively.

We are interested in identifying the conditions under which populations transition from a steady state to a chaotic state. To model the steady state, we employ the gamma distribution, which is widely acknowledged as an effective approximation for stationary distributions in various ecological contexts \cite{Dennis1984, Pielou1975, Engen1978}. This approach has proven especially useful in studies of species such as the \textit{Tribolium} beetle, where laboratory populations exhibit fluctuations around a mean equilibrium size due to environmental variability \cite{Dennis1984, Peters1989, DennisCostantino1988, Costantino1981}. These results not only highlight the gamma distribution’s practical utility for stochastic modeling but also provide a foundation for exploring the transition mechanisms from steady to chaotic dynamics.

We first establish the mathematical relations between the intrinsic growth rate \( r \) in \eqref{eq:logistic3} and \eqref{ricker1} and the variance of \( \epsilon_t \) along with the parameters of the gamma distribution at equilibrium, which allow us to analyze the range of $r$ for various \( k \) and \( \epsilon_t \).  The relations reveal that, both \eqref{eq:logistic3} and \eqref{ricker1} have two branches of the intrinsic growth rate, \( r_+ \) and \( r_- \), representing alternative stable states corresponding to higher and lower growth rates. The findings indicate that, in the logistic model, both of \( r_+ \) and \( r_- \) are less than $3$, indicating there is no transition between a steady state and a chaotic state. On the other hand, in the Ricker model, for all \( \epsilon_t \), there is a transition from a steady state and a chaotic state when the shape parameter $k$ of the gamma distribution is small.
 
Building on the work in \cite{wang2025}, this study continues to study stochastic equations with a focus on transitions from steady states to chaos. We begin by exploring the interplay between deterministic chaos and stochastic dynamics and introduce a novel method for identifying the shift from steady-state behavior to chaotic regimes. This approach leverages classical chaos theory results in conjunction with the gamma distribution to model the underlying population dynamics. Our analysis reveals a striking difference between the stochastic logistic and Ricker models in terms of how they transition from steady states to chaos. Recognizing these transitions is crucial, with significant implications across fields such as ecology, biology, climate science, and economics.

\section{Stochastic equations converge to deterministic equations}

\subsection{The Gamma Distribution}

The gamma distribution is defined on the interval $(0, \infty)$ and is characterized by two positive parameters: the shape parameter $k$ and the scale parameter $\theta$. Its probability density function (PDF) is:

\begin{equation} \label{eq:gamma_pdf}
f(x; k, \theta) = \frac{x^{k - 1} e^{-x / \theta}}{\Gamma(k) \theta^k}, \quad \text{for } x > 0,
\end{equation}
where $\Gamma(\cdot)$ is the gamma function. For positive integer $n$, $\Gamma(n)=(n-1)!$.  Within the gamma distribution, the parameter $k$, often referred to as the shape parameter, influences the form and skewness of the distribution, while the scale parameter $\theta$ governs its dispersion. This distribution finds applications in various fields, including population dynamics, event waiting times, and biological processes \cite{wikipedia_gamma}.

Many studies have demonstrated that the gamma distribution accurately approximates these stationary distributions in diverse ecological settings \cite{Dennis1984,Peters1989,DennisCostantino1988,Costantino1981}. In our analysis, we assume that the population size $X_t$ in \eqref{eq:logistic3} and \eqref{ricker1} at equilibrium follows a gamma distribution with parameters $k$ and $\theta$:

\begin{equation} \label{eq:x_gamma_distribution}
X_t \sim \text{Gamma}(k, \theta).
\end{equation}
The mean and variance of the gamma distribution are:
\begin{align}
\mu &= E[X_t] = k \theta, \label{eq:gamma_mean} \\
\sigma^2 &= \text{Var}[X_t] = k \theta^2. \label{eq:gamma_variance}
\end{align}

\subsection{Mathematical relation between deterministic and stochastic equations}
The following relation between stochastic and deterministic equations indicates that the stochastic equations converge to their corresponding deterministic equations when the variance of $X_t$ is sufficient small. Here the requirement of the gamma distribution on $X_t$ can be relaxed.

\begin{theorem}
Let $X_t$ be a random variable representing the population at time $t$. Let $y_t=E[X_t]$. Assume that the steady state of stochastic discrete model \eqref{eq:logistic3}, $X_t$, follows the gamma distribution \eqref{eq:x_gamma_distribution}.  Then stochastic logistic equation \eqref{eq:logistic3} implies that 
\begin{equation}
    y_{t+1} =  ry_t(1-y_t) - r\text{Var}(X_t)
\label{eq:relation1}
\end{equation}
And stochastic Ricker equation \eqref{ricker1} implies that 
\begin{equation}
    y_{t+1} =  y_t\,\exp\!\Bigl(r(1-y_t)\Bigr) + O\Big(\operatorname{Var}(X_t)\Bigr).
\label{eq:relation2}
\end{equation}
\label{theorem0}
\end{theorem}

{\bf Proof of Theorem \ref{theorem0}} Let's start with stochastic equation \eqref{eq:logistic3}. If $X_t$ satisfies \eqref{eq:logistic3}, taking the expectation at the both sides implies that
\begin{align}
    y_{t+1} &= E[rX_t(1-X_t) \epsilon_t] \nonumber \\
    &= rE[X_t] - rE[X_t^2] \nonumber\\
    &= ry_t - r(y_t^2 + \text{Var}(X_t)) \nonumber\\
    &= ry_t - ry_t^2 - r\text{Var}(X_t) \nonumber\\
    &= ry_t(1-y_t) - r\text{Var}(X_t)
\end{align}
which is \ref{eq:relation1}. Apparently, the gamma distribution assumption is not used in the proof for \ref{eq:relation1}.   Now if \(X_t\) satisfies stochastic Ricker equation \eqref{eq:logistic3}, then 
\begin{equation}
    X_{t+1} = X_t \exp\!\Bigl(r(1-X_t)\Bigr) \epsilon_t,
\end{equation}
and that \(X_t\) is a Gamma-distributed random variable with parameters \(k\) and \(\theta\), so that
\[
E[X_t] = y_t = k\,\theta \quad \text{and} \quad \operatorname{Var}(X_t) = k\,\theta^2.
\]
We write
\begin{equation}
    X_t = y_t + \beta_t,\quad \text{with } E[\beta_t]=0.
\end{equation}

Then,
\begin{equation}
    X_{t+1} = X_t \exp\!\Bigl(r(1-X_t)\Bigr)  \epsilon_t
    = (y_t+\beta_t)\exp\!\Bigl(r\Bigl(1-(y_t+\beta_t)\Bigr)\Bigr)  \epsilon_t.
\end{equation}
Rewrite the exponential as
\begin{equation}
    \exp\!\Bigl(r\Bigl(1-(y_t+\beta_t)\Bigr)\Bigr)
    = \exp\!\Bigl(r(1-y_t)\Bigr)\exp\!\Bigl(-r\beta_t\Bigr).
\end{equation}

Expand \(\exp(-r\beta_t)\) about \(\beta_t=0\) up to third order:
\begin{equation}
\label{eq:expansion}
    \exp(-r\beta_t) = 1 - r\beta_t + \frac{1}{2}(r\beta_t)^2 - \frac{1}{6}(r\beta_t)^3 + R_3(\beta_t),
\end{equation}
where the remainder \(R_3(\beta_t)\) is of order \(\beta_t^4\).

Thus, we have
\begin{align}
    X_{t+1} &= (y_t+\beta_t)\,\exp\!\Bigl(r(1-y_t)\Bigr)
    \Bigl[ 1 - r\beta_t + \frac{1}{2}(r\beta_t)^2 - \frac{1}{6}(r\beta_t)^3 + R_3(\beta_t) \Bigr]  \epsilon_t.
\end{align}
Taking the expectation of both sides and using the linearity of expectation, we obtain
\begin{align}
    E[X_{t+1}] &= \exp\!\Bigl(r(1-y_t)\Bigr) \Biggl\{ y_t \nonumber\\[1mm]
    &\quad + \Bigl[- y_t\,E[\beta_t] + E[\beta_t]\Bigr] \nonumber\\[1mm]
    &\quad + \Bigl[\frac{1}{2} y_t\,E[(r\beta_t)^2] - E[(r\beta_t)^2]\Bigr] \nonumber\\[1mm]
    &\quad + \Bigl[-\frac{1}{6}y_t\,E[(r\beta_t)^3] + \frac{1}{2}E[(r\beta_t)^3]\Bigr] \nonumber\\[1mm]
    &\quad + E\Bigl[(y_t+\beta_t)R_3(\beta_t)\Bigr] \Biggr\}.
\end{align}
Since \(E[\beta_t]=0\), the first-order term cancels. Moreover, noting that, 
\[
E[(r\beta_t)^2] = r^2 E[\beta_t^2] = r^2\,\operatorname{Var}(X_t),
\]
and that the third-order term is \(O((\operatorname{Var}(X_t))^{3/2})\) from \eqref{estimation1}, and from \eqref{estimation2}, 
\[
E\Bigl[(y_t+\beta_t)R_3(\beta_t)\Bigr] = O\Bigl((\operatorname{Var}(X_t))^2\Bigr),
\]
we obtain the final approximate expression:
\begin{align}
    y_{t+1} &= E[X_{t+1}] = y_t \exp\!\Bigl(r(1-y_t)\Bigr)\\ 
    &\quad + \exp\!\Bigl(r(1-y_t)\Bigr)\left(\frac{r^2 y_t}{2} - r\right) \operatorname{Var}(X_t) \nonumber\\
    &\quad + O\Bigl((\operatorname{Var}(X_t))^{3/2}\Bigr)+ O\Bigl((\operatorname{Var}(X_t))^2\Bigr) \nonumber.
\end{align}
which is \eqref{eq:relation2}.

\section{Transitions from steady states to chaos}
Theorem \ref{theorem0} indicates that when the variance of $X_t$ is sufficient small, the stochastic equations converge to the corresponding deterministic equation for their expectations with the same intrinsic growth rate \( r \).  Here we explore the impact of the variance of \( \epsilon_t \) and $k$ on transitions from a steady state to a chaotic state  of \eqref{eq:logistic3} and \eqref{ricker1}. 

Due to the nonlinearity in equations (\ref{eq:logistic3}) and (\ref{ricker1}), $X_{t+1}$ may not necessarily follow the gamma distribution. However, it is reasonable to assume that the population at equilibrium from specific time $t$ to $t+1$ maintains the same expectation and variance. Therefore, at equilibrium, we assume that, for a specific time $t$,  
\begin{equation} \label{eq:conditions}
E[X_{t+1}] = E[X_t], \,\,\,  Var[X_{t+1}] = Var[X_t] 
\end{equation}
which allows us to derive explicit mathematical relations for $r$ in terms of $k$ and $\theta$.

\subsection{Results on logistic equation}
Here we state the mathematical relation among the intrinsic growth rate \( r \) in \eqref{eq:logistic3}, the variance of \( \epsilon_t \) along with the parameters of the gamma distribution on equilibrium states. This simplifies the result in \cite{wang2025} and its proof is in the Appendix \ref{logisticsection}.   In particular, our discussion focuses on the impact of $\text{Var}[\epsilon_t] $ on $r$ for the study of transition from steady states to chaos.    
\begin{theorem}
Let the steady state of stochastic discrete model \eqref{eq:logistic3}, $X_t$, follow the gamma distribution \eqref{eq:x_gamma_distribution}, and its expectation and variance remain constant from $t$ to $t+1$ \eqref{eq:conditions}. Then \( r \) in \eqref{eq:logistic3} is determined in terms of $k,\text{Var}[\epsilon_t] $ as 
\begin{align}
r_{\pm} = \dfrac{ (2k + 4) \pm (k + 1) \sqrt{ \dfrac{1 - Var[\epsilon_t](k + 2)}{Var[\epsilon_t] + 1} } }{ k + 3 }.
\label{eq:multiplied_q31}
\end{align}
\label{theorm1}
\end{theorem}
When the variance of $\epsilon_t$  is zero, (\ref{eq:multiplied_q31}) has one trivial solution $r=1$ and 
\begin{equation}
r = \frac{ 3k + 5 }{ k + 3 } 
\label{eq:expression_r}
\end{equation}
From equation (\ref{eq:expression_r}), it follows that 
$$
\frac{5}{3}<r<3. 
$$
It is well known that the population dynamics of the deterministic logistic model (\ref{eq:dlogistic3}), such as stable growth, periodic oscillations, chaos, or extinction, are dependent on the value of \( r \). If \( 1 < r < 3 \), the population grows and eventually reaches a nonzero steady state\cite{may1976simple,strogatz2018nonlinear,hilborn2000chaos}. For the stochastic logistic equation (\ref{eq:logistic3}), the feasible range of $r$, $ (\frac{5}{3}, 3) \subset (1,3) $, at equilibrium implies that there is no transition from a steady state to chaotic state. 


\subsubsection{Impact of $k$ and $\text{Var}(\epsilon_t)$}

\begin{figure}[h!]
    \centering
    \begin{subfigure}[b]{0.45\textwidth}
        \centering
        \includegraphics[width=\textwidth]{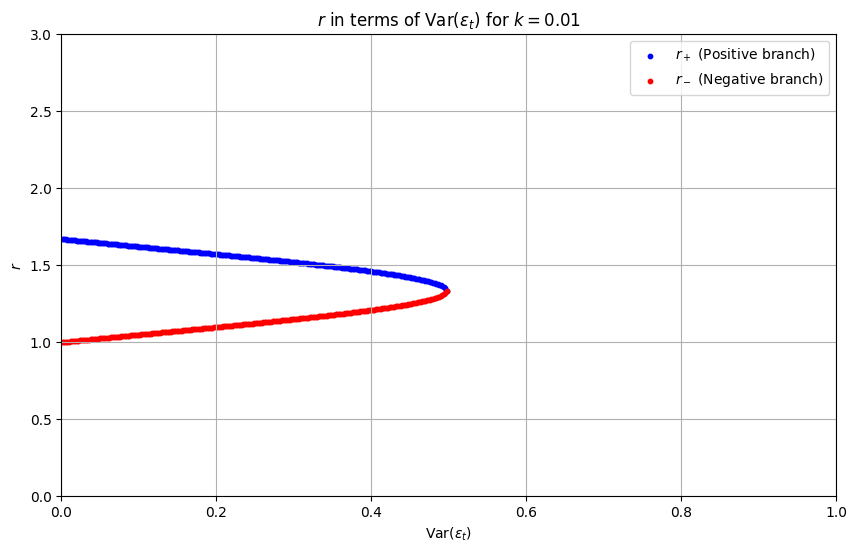}
        \label{fig:sub1}
    \end{subfigure}
    \hfill
    \begin{subfigure}[b]{0.45\textwidth}
        \centering
        \includegraphics[width=\textwidth]{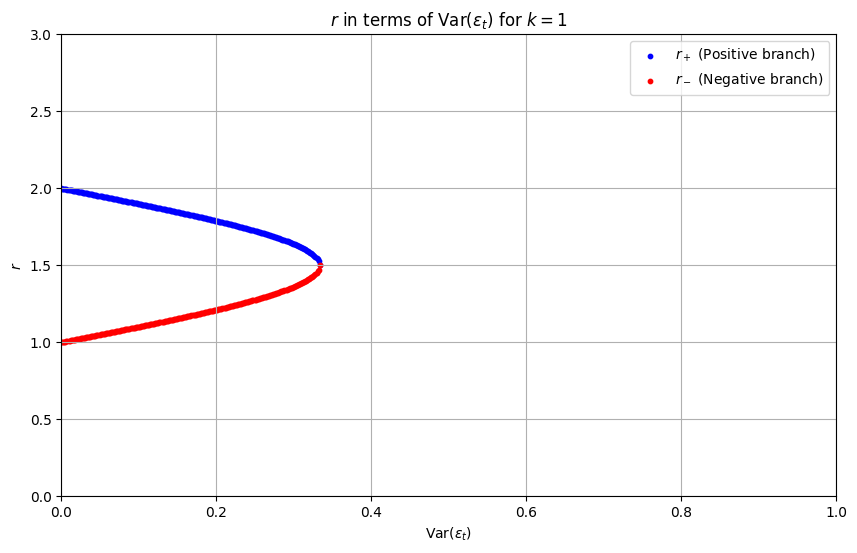}
        \label{fig:sub2}
    \end{subfigure}
    
    \vskip\baselineskip
    
    \begin{subfigure}[b]{0.45\textwidth}
        \centering
        \includegraphics[width=\textwidth]{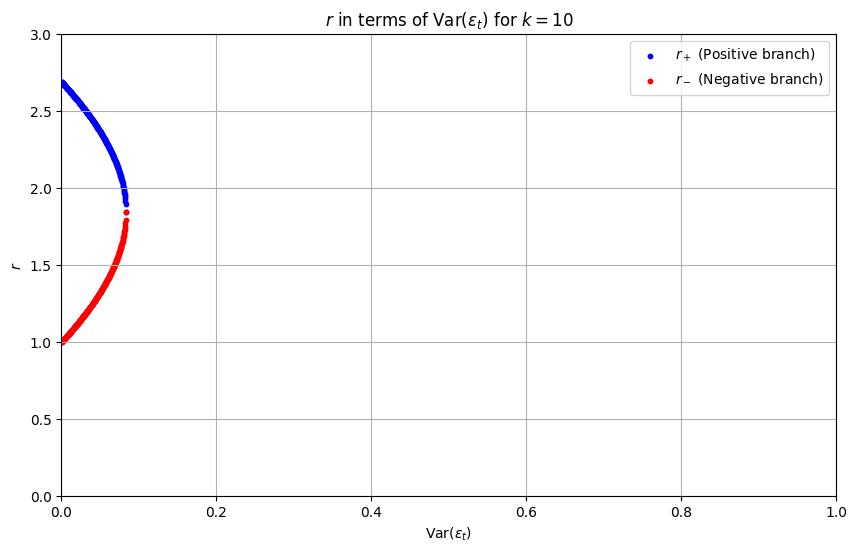}
        \label{fig:sub3}
    \end{subfigure}
    \hfill
    \begin{subfigure}[b]{0.45\textwidth}
        \centering
        \includegraphics[width=\textwidth]{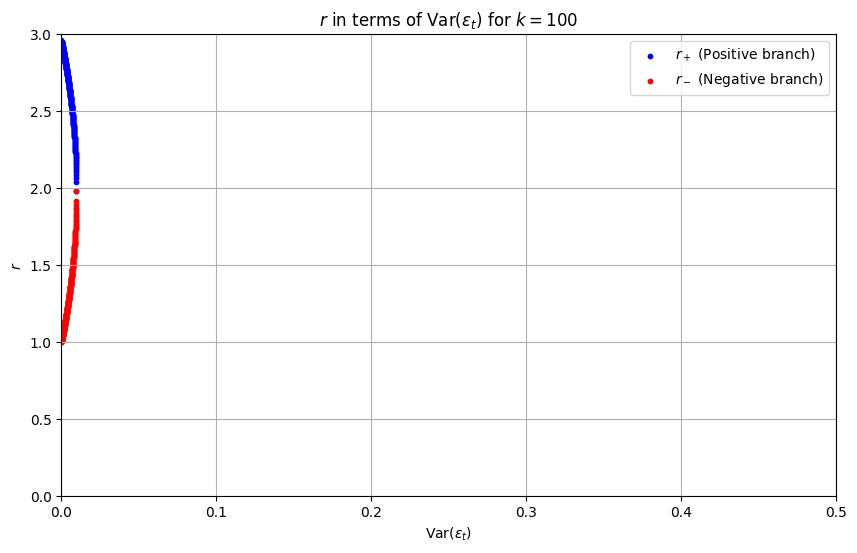}
        \label{fig:sub4}
    \end{subfigure}
    
    \caption{Plots of $r$ in term of $Var[\epsilon_t]$}
    \label{fig:main}
\end{figure}

By Theorem~\ref{theorm1}, the parameters at equilibrium must satisfy
\[
0 < \mathrm{Var}[\epsilon_t] \leq \frac{1}{k+2}
\quad \text{and} \quad
0 \leq \mathrm{Var}[\epsilon_t] \leq 0.5.
\]
To illustrate how these constraints affect \(r\), we present simulations in Figure~\ref{fig:main}.
The four panels show the relationship between the growth rate \(r\) and 
\(\mathrm{Var}(\epsilon_t)\) for different values of the parameter \(k\).
  
First, we observe that all values of \( r \) in Figure~\ref{fig:main} fall within the range \((1,3)\). 
For the deterministic logistic model~\eqref{eq:dlogistic3}, populations grow and stabilize in precisely this interval, 
where the growth rate balances the carrying capacity, enabling the population to reach a stable state 
\cite{may1976simple,strogatz2018nonlinear,hilborn2000chaos}. 
As a result, there is no transition from a steady state to a chaotic state in~\eqref{eq:logistic3}.

Larger values of $k$ (more symmetrical and regular population distribution) extends the range of $r$ as indicated in Figure \ref{fig:main}. For larger values of \( k \), such as \( \text{Var}(\epsilon_t) = 10 \) and \(k = 100 \), the feasible range of \( r \) is extensive. Populations in this setting are expected to exhibit a wide range of equilibrium growth rates, with both \( r_+ \) and \( r_- \) branches present.  Larger values of \(k\) (indicating a more symmetrical, regular population distribution)
extend the range of \(r\), as shown in Figure~\ref{fig:main}. 
For instance, when \(\mathrm{Var}(\epsilon_t) = 10\) and \(k = 100\), the feasible range of \(r\) becomes very broad.
In such a scenario, populations can exhibit a wide range of equilibrium growth rates, with both \(r_+\) and \(r_-\) 
branches present. A high \(k\) value coupled with low variance indicates that the population distribution 
is less skewed and more resilient, allowing for a stable equilibrium across an extensive range of \(r\).
Populations in this context resemble species inhabiting stable environments with minimal environmental variability, 
thereby supporting predictable population dynamics around the equilibrium.

As \( k \) decreases (e.g., \( k = 0.01 \) or \( k = 1 \)), the range of \( r \) values that produce stable solutions 
becomes narrower. This suggests that populations with a more skewed distribution have a limited ability 
to attain stable equilibria.

\subsection{Results on Ricker Equation}
Here we state the mathematical relation among the intrinsic growth rate \( r \) in \eqref{ricker1} and the variance of \( \epsilon_t \) along with the parameters of the gamma distribution on equilibrium states. This basically is the result in \cite{wang2025}. Our discussion focuses on the impact of $\text{Var}[\epsilon_t] $ on $r$ for the study of transition from steady states to chaos.  
\begin{theorem}
Let the steady state of stochastic discrete model \eqref{ricker1},  $X_t$, follow the gamma distribution \eqref{eq:x_gamma_distribution}, and its expectation and variance remain constant from $t$ to $t+1$ \eqref{eq:conditions}. Then \( r \) \eqref{ricker1} depends on $k, \text{Var}[\epsilon_t]$ with the following the relation
\begin{equation}
2 e^{\frac{r}{k + 1}} - \left ((1+Var[\epsilon_t])e^{2r}\right )^{\frac{1}{k + 2}} = 1
\label{final_equation1}
\end{equation}
\label{theorm2}
\end{theorem}

\subsubsection{$Var[\epsilon_t]=0$}
\begin{figure}[h!]
    \centering
    \includegraphics[width=0.6\textwidth]{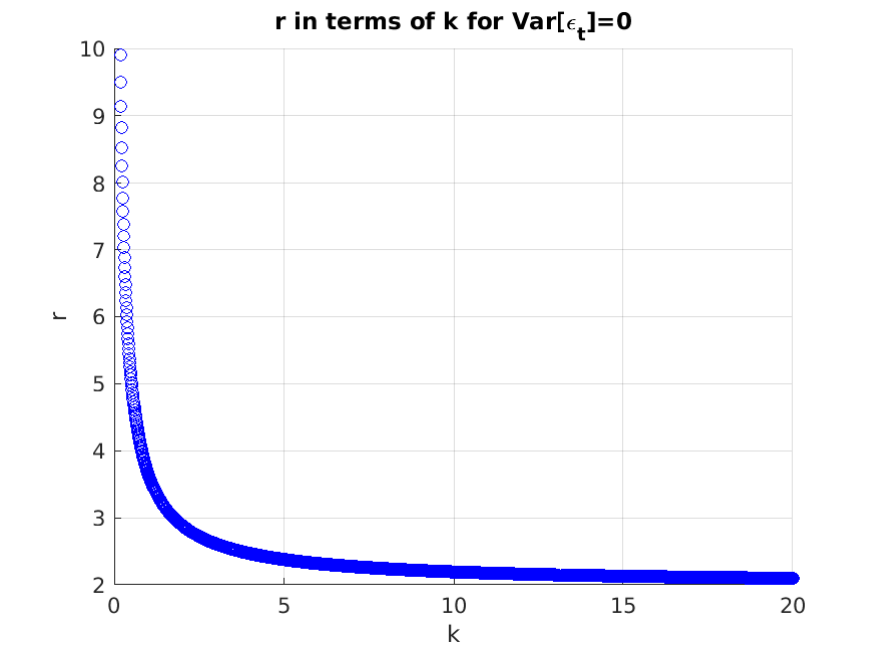}
    \caption{Positive $r$ in terms of $k$ for $Var[\epsilon_t]=0$ ($r=0$ is a solution.) }
    \label{plot11}
\end{figure}
When the variance of $\epsilon_t$  is zero, equation \eqref{final_equation1} becomes
\begin{equation}
2 e^{\frac{r}{k + 1}} - \left (e^{2r}\right )^{\frac{1}{k + 2}} = 1
\label{eq:multiplied_nr5}
\end{equation}
Note that \(r = 0\) is always a trivial solution of \eqref{eq:multiplied_nr5}. 
We illustrate the positive solution \(r\) of \eqref{eq:multiplied_nr5} as a function of \(k\) 
in Figure~\ref{plot11}. 
It is evident that \(r\) decreases with respect to \(k\), while remaining above 2. 

For the deterministic Ricker model \eqref{dricker1}, the population dynamics depends intricately on the value of \(r\) \cite{may1976simple,kot2001elements}.  When \(r > 2\), the system \eqref{dricker1} becomes unstable, leading to cycles of period \(2^n\), 
and eventually transitions into a chaotic regime \cite{may1976simple,strogatz2018nonlinear,hilborn2000chaos}. 
Hence, for \(\mathrm{Var}[\epsilon_t] = 0\), the stochastic discrete model \eqref{ricker1} may admit chaotic behavior.

\subsubsection{Impact of $k$ and $Var[\epsilon_t$]}
\begin{figure}[h!]
    \centering
    \begin{subfigure}[b]{0.45\textwidth}
        \centering
        \includegraphics[width=\textwidth]{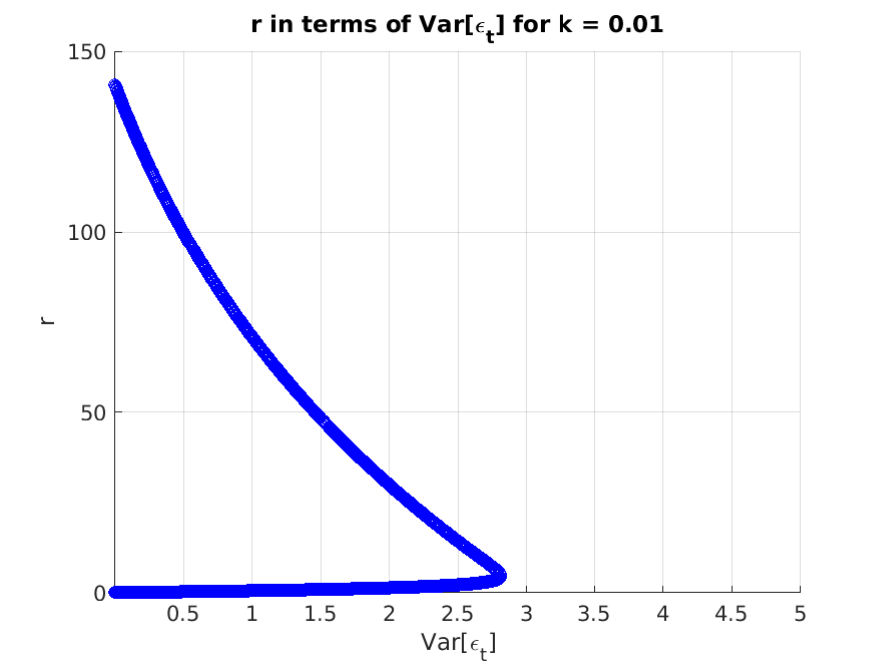}
        \label{fig:sub1}
    \end{subfigure}
    \hfill
    \begin{subfigure}[b]{0.45\textwidth}
        \centering
        \includegraphics[width=\textwidth]{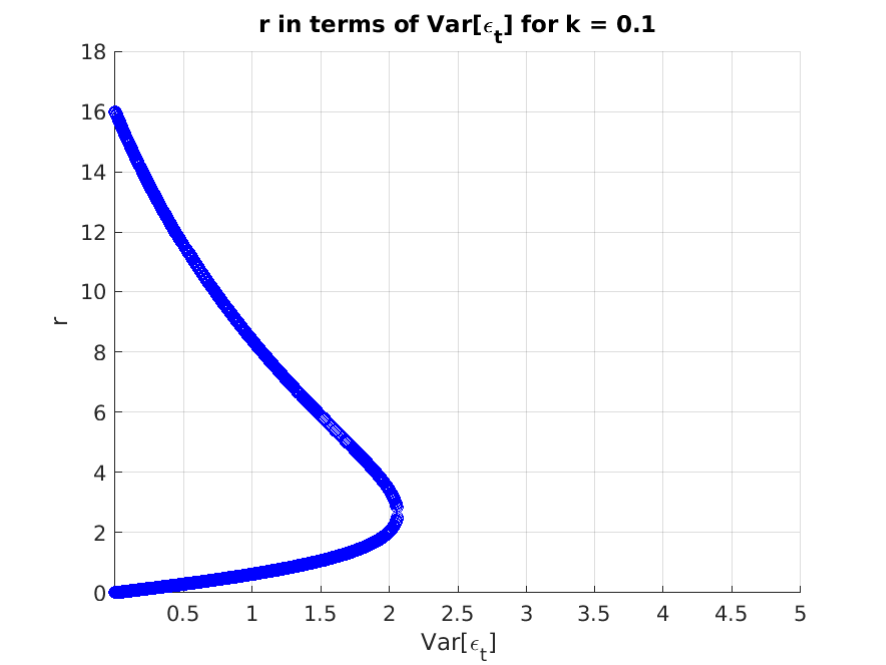}
        \label{fig:sub2}
    \end{subfigure}
    
    \vskip\baselineskip
    
    \begin{subfigure}[b]{0.45\textwidth}
        \centering
        \includegraphics[width=\textwidth]{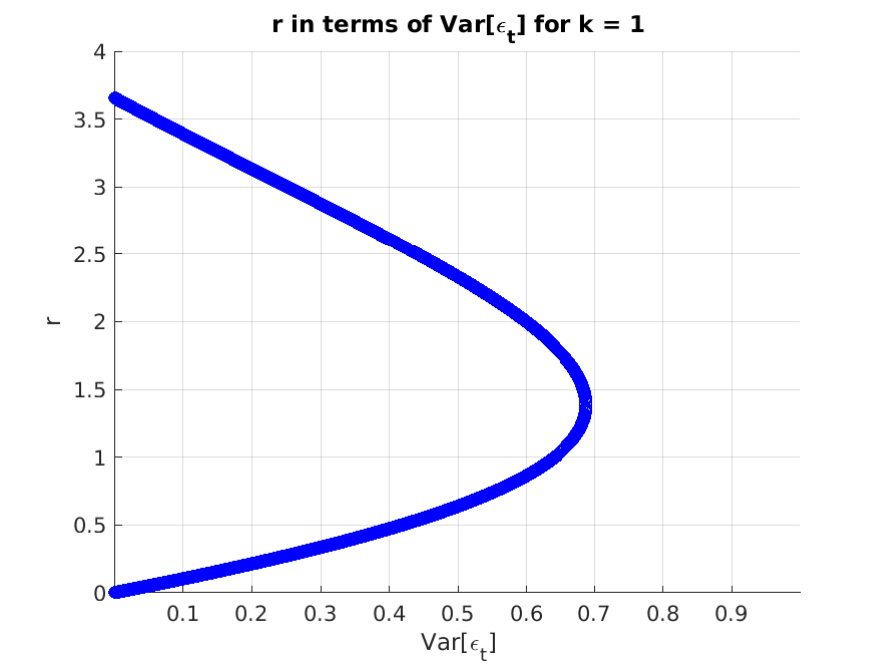}
        \label{fig:sub3}
    \end{subfigure}
    \hfill
    \begin{subfigure}[b]{0.45\textwidth}
        \centering
        \includegraphics[width=\textwidth]{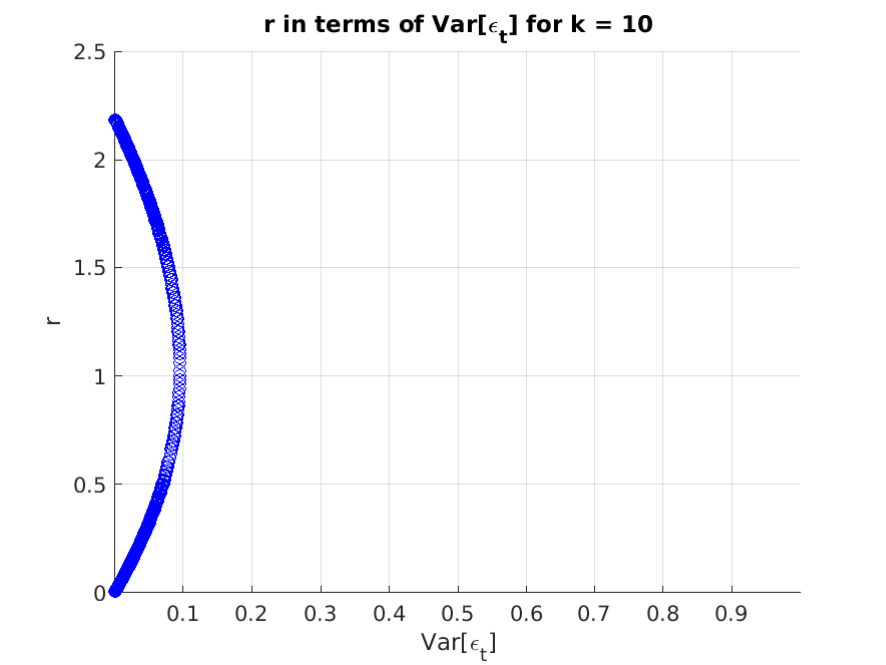}
        \label{fig:sub4}
    \end{subfigure}
        
    \caption{Plots of $r$ in terms of $Var[\epsilon_t]$ for Ricker equation}
    \label{fig:main1}
\end{figure}

In this section, we employ numerical simulations to explore the relationship between \(r\) and \(k\) in 
\eqref{final_equation1} for the stochastic Ricker equation~\eqref{ricker1}. 
Our aim is to reveal how variations in \(k\) and \(\mathrm{Var}[\epsilon_t]\) influence \(r\). 
The results, presented in Figure~\ref{fig:main1}, comprise four panels illustrating the relationship 
between the growth rate \(r\) and \(\mathrm{Var}[\epsilon_t]\) for different values of \(k\).
    
For \(\mathrm{Var}[\epsilon_t] > 0\), there are two branches of \(r\). The upper bounds of 
\(\mathrm{Var}[\epsilon_t]\) are determined by \(k\), and smaller values of \(k\) extend these upper bounds, 
as indicated in Figure~\ref{fig:main1}.  While the lower branch of \( r \) is less than $2$ at equilibrium, the upper branch of $r$ of the stochastic Ricker equation \eqref{ricker1}, may exceed 2 for $k$ is small, and enter a chaotic region, characterized by unpredictable and irregular fluctuations.  Smaller values of $k$ correspond to highly skewed distributions. Populations with low $k$ values are more susceptible to environmental fluctuations.

Larger values of $k$ (more symmetrical and regular population distribution) shrink the range of $r$ as indicated in Figure \ref{fig:main1}, which can reduce the possibility to enter the transition region from s steady state to a chaotic state. For larger values of \( k \), such as \(k = 100 \), the feasible range of \( r \) is limited to $ (0, 2.5)$. Populations in this setting are expected to exhibit a wide range of equilibrium growth rates, with both \( r_+ \) and \( r_- \) branches present.  Populations with this setup are analogous to species in stable environments that experience minimal environmental variability, allowing them to maintain predictable population dynamics around the equilibrium.

As \( k \) decreases (e.g., \( k = 0.01 \) and \( k = 1 \)), the range of \( r \) values for stable solutions extends and enters the chaotic region.  This indicates that, even with small $\text{Var}[\epsilon_t]$, highly skewed distributions population could have a chaotic state.

\section{Discussion}
Chaos represents a fundamental change in system dynamics, often involving complex processes that are both theoretically and practically fascinating. The transition between a steady state and chaos is particularly interesting because it reveals the complex, nonlinear nature of many systems. It highlights how small changes can lead to unpredictable, large-scale dynamics. Understanding this transition is crucial for managing ecosystems, predicting climate patterns, and controlling chaotic behaviors in technology and biology \cite{earlywarning2023}.

In this paper, we first demonstrate that the interplay between deterministic and stochastic equations, and then present  a novel approach for identifying such a transition from a steady state to chaos based on the classical chaos results and the gamma distribution for the underlying population. Our findings using the gamma distribution in this paper reveals that, both \eqref{eq:logistic3} and \eqref{ricker1} have two branches of the intrinsic growth rate, \( r_+ \) and \( r_- \), representing alternative stable states corresponding to higher and lower growth rates. However, there is a striking difference in between the logistic model and Ricker model: there is no transition between a steady state and a chaotic state for the logistic model, but there is a transition from a steady state and a chaotic state when the shape parameter $k$ of the gamma distribution is small. 

The striking difference for the transition from steady states to chaos between of the stochastic logistic and Ricker model even exhibits in deterministic equations. In fact,  the Ricker model more readily exhibits chaotic behavior compared to the logistic model, primarily due to its nonlinear density dependence and exponential growth suppression. The logistic model, given by \(x_{n+1} = r x_n (1 - x_n)\), assumes linear density dependence, where population growth slows as the population approaches the carrying capacity \cite{may1976simple}. In this model, chaos emerges at higher growth rates (\(r > 3.57\)) following a period-doubling route to chaos \cite{feigenbaum1980quantitative}. In contrast, the Ricker model, defined as \(x_{n+1} = x_n e^{r(1 - x_n)}\), incorporates exponential density dependence, leading to sharper suppression of population growth at higher densities \cite{ricker1954stock}. This sharper feedback allows the Ricker model to exhibit chaotic dynamics at much lower growth rates (\(r \sim 2.5\)) compared to the logistic model. The exponential term in the Ricker model results in faster and more intense oscillations, making chaotic behavior more prevalent in systems such as fish populations or insect dynamics \cite{may1974biological}. Consequently, while the logistic model provides a simpler framework for understanding chaos, the Ricker model is often better suited for studying real-world populations where overshooting and collapses are common.

Automated methods for chaos detection in noisy empirical data have also been developed, addressing the limitations of classical approaches that are sensitive to noise \cite{arxiv2019chaos}. Furthermore, combining deterministic and stochastic dynamics provides a framework for understanding extreme events in chaotic systems, such as those observed in climate variability and ecological dynamics \cite{chaos2017extremes}. Advanced methods, such as deep learning, are now being employed to differentiate chaotic and stochastic behavior in biological systems, including lineage trees, offering robust tools for chaos detection in noisy and high-dimensional data \cite{physrev2022lineagetrees}. 

The method with the gamma distribution in this paper provide complementary approaches to detecting chaos in biological systems. Based on the theoretical relationships between stochastic models and statistical distributions, future research could leverage laboratory and real-world ecological population data to analyze parameter-dependent transitions to chaos and differentiate chaos from noise or periodicity.  This process would involve fitting the stochastic models to empirical data to estimate their parameters, as well as deriving the corresponding distribution parameters. By applying this methods along with other relevant chaos detection methods, researchers can uncover deterministic chaos, informing fields such as population ecology, neuroscience, and cardiology. Future research should also focus on improving these techniques' robustness to noise and their applicability to high-dimensional systems.

\section{Appendix}
\subsection{Expectation and Variance Conditions of logistic Equation}\label{logisticsection}
Here we derive the equilibrium conditions for the mean and variance of~\eqref{eq:logistic3}, assuming that \( X_t \) is gamma-distributed at equilibrium. That is, at equilibrium, \( X_t \sim \mathrm{Gamma}(k, \theta) \), where \( k \) is the shape parameter and \( \theta \) is the scale parameter:

\begin{equation}
X_t \sim \text{Gamma}(k, \theta)
\label{gamma123}
\end{equation}
Because of nonlinearity, $X_{t+1}$ may not follow a gamma distribution, it is reasonable to assume that the mean and variance are the same: 
\begin{equation}
E[X_{t+1}] = E[X_t] \quad \text{and} \quad \text{Var}[X_{t+1}] = \text{Var}[X_t]
\label{mean_identity-123}
\end{equation}
It is easy to compute that 
\begin{align}
E[X_t] &= \mu = k \theta, \label{eq:mean_mu} \\
E[X_t^2] &= k (k + 1) \theta^2. \label{eq:second_moment}
\end{align}
\begin{align}
E[X_t^n] &= \theta^n \frac{\Gamma(k+n )}{\Gamma(k)}, \label{eq:gamma_moment} \\
\end{align}
Since $\epsilon_t$ is independent, we have  
\begin{equation} \label{eq:expected_value}
E[X_{t+1}] = E\left[ r X_t \left(1 - X_t \right )  \epsilon_t \right]=E\left[ r X_t \left(1 - X_t \right ) \right] E[\epsilon_t].
\end{equation}
Since $E[\epsilon_t]=1$, and at equilibrium, ($E[X_{t+1}] = E[X_t] = \mu$), we have
\begin{equation} \label{eq:expanded_equation}
\mu = r \left( \mu - E[X_t^2] \right).
\end{equation}
Subtract $r \mu$ from both sides and simplify it:
\begin{equation} \label{eq:simplified_equation}
\mu (1 - r) = - r E[X_t^2].
\end{equation}
Substituting $\mu = k \theta$ and $E[X_t^2]=\sigma^2 + \mu^2$ into the equilibrium equation \eqref{eq:simplified_equation}:
\begin{equation} \label{eq:equilibrium_in_k_theta}
k \theta (1 - r) = - r \left( k \theta^2 + k^2 \theta^2 \right).
\end{equation}
Divide both sides by $k \theta$ (since $k, \theta > 0$):

\begin{equation} \label{eq:equilibrium_divided}
(1 - r) = - r \theta (1 + k).
\end{equation}
Solve for $\theta$:

\begin{equation} \label{eq:theta_solution}
\theta = \frac{r - 1}{r (1 + k)}.
\end{equation}
Therefore, this suggests that for $r > 1$, this equilibrium may exist under the gamma distribution assumption.

In addition to the mean, we can compare the variance condition at equilibrium $$ \text{Var}[X_{t+1}] = \text{Var}[X_t].$$ This provides additional conditions that can help determine the gamma distribution parameters. To compute \( \text{Var}[X_{t+1}] \), we need \( E[X_{t+1}] \) and \( E[X_{t+1}^2] \):

\[
\text{Var}[X_{t+1}] = E[X_{t+1}^2] - \left( E[X_{t+1}] \right)^2.
\]
We already have \( E[X_{t+1}] = \mu = k \theta \). To compute \( E[X_{t+1}^2] \), we expand \( X_{t+1}^2 \):

\begin{align*}
X_{t+1}^2 &= \left( r X_t - r X_t^2  \right)^2 \epsilon_t^2 \\
&= \left (r^2 X_t^2 - 2 r^2 X_t^3 + r^2 X_t^4 \right )\epsilon_t^2
\end{align*}
We take expectations:
\begin{equation} \label{eq:expectation_xt1_squared}
E[X_{t+1}^2] = \left ( r^2 E[X_t^2] - 2 r^2 E[X_t^3 ] + r^2 E[X_t^4] \right )E[\epsilon_t^2].
\end{equation}
At equilibrium, \( \text{Var}[X_{t+1}] = \text{Var}[X_t] = k \theta^2 \). Therefore,

\begin{equation} \label{eq:variance_equilibrium}
E[X_{t+1}^2] - (k \theta)^2 = k \theta^2
\end{equation}
and
\begin{equation} \label{eq:variance_equilibrium_simplified}
E[X_{t+1}^2] = k \theta^2 (1 + k).
\end{equation}
We can cancel \( \theta^2 \) from both sides as \( \theta > 0 \):

\begin{equation} \label{eq:variance_equilibrium_theta_cancelled}
\frac{E[X_{t+1}^2]}{\theta^2} = k (1 + k).
\end{equation}
Compute the required terms:

\begin{align*}
E[X_t^2] &= \theta^2 k (k + 1), \\
E[X_t^{3}] &= \theta^{3} \frac{\Gamma(k + 3)}{\Gamma(k)} = \theta^{3} (k + 2)(k + 1)k, \\
E[X_t^{4}] &= \theta^{4} \frac{\Gamma(k + 4)}{\Gamma(k)} = \theta^{4} (k + 3)(k + 2)(k + 1)k.
\end{align*}
With cancellation of $\theta^2$,  \eqref{eq:variance_equilibrium_theta_cancelled} and \eqref{eq:expectation_xt1_squared} give
\begin{align*}
& \Big( r^2 k (k + 1) - 2 r^2 \theta (k + 2)(k + 1)k + r^2 \theta^{2}  (k + 3)(k + 2)(k + 1)k \Big ) E[\epsilon_t^2] = k (1 + k).
\end{align*}
Dividing by $E[\epsilon_t^2]$ and subtracting $k (k + 1)$ at both sides give  
\begin{align*}
&  (r^2-1) k (k + 1) - 2 r^2 \theta (k + 2)(k + 1)k + r^2 \theta^{2}  (k + 3)(k + 2)(k + 1)k = \frac{k (1 + k)}{E[\epsilon_t^2]}-k(k+1).
\end{align*}
Multiplying by $-1$ at both sides produces 
\begin{align*}
&  (1-r^2) k (k + 1) + 2 r^2 \theta (k + 2)(k + 1)k - r^2 \theta^{2}  (k + 3)(k + 2)(k + 1)k = \frac{E[\epsilon_t^2]-1}{E[\epsilon_t^2]}k(k+1).
\end{align*}
Since $Var[\epsilon_t]=E[\epsilon_t^2]-E[\epsilon_t]^2=E[\epsilon_t^2]-1$, 
\begin{align*}
&  (1-r^2) k (k + 1) + 2 r^2 \theta (k + 2)(k + 1)k - r^2 \theta^{2}  (k + 3)(k + 2)(k + 1)k =\frac{Var[\epsilon_t]}{Var[\epsilon_t]+1}k(k+1).
\end{align*}
%
%
%
%
%
Note that $\theta = \frac{(r - 1)}{r (1 + k)}$, therefore,
\begin{equation}
(1 - r^2) k (k + 1) + 2 r (r - 1) k (k + 2) - (r - 1)^2 k \cdot \dfrac{ (k + 2)(k + 3) }{ k + 1 } = \frac{Var[\epsilon_t]}{Var[\epsilon_t]+1}k(k+1)
\label{eq:simplified_main_n1}
\end{equation}
Dividing both sides by \( k (r - 1) \), we have

\begin{equation}
\dfrac{ (1 - r^2) k (k + 1) }{ k (r - 1) } + \dfrac{ 2 r (r - 1) k (k + 2) }{ k (r - 1) } - \dfrac{ (r - 1)^2 k \cdot \dfrac{ (k + 2)(k + 3) }{ k + 1 } }{ k (r - 1) } = \frac{Var[\epsilon_t]}{Var[\epsilon_t]+1}\frac{k+1}{r - 1}
\label{eq:divided_n1}
\end{equation}
After simplification, equation \eqref{eq:divided_n1} becomes:
\begin{equation}
- (r + 1)(k + 1) + 2 r (k + 2) - \dfrac{ (r - 1) (k + 2)(k + 3) }{ k + 1 } = \frac{Var[\epsilon_t]}{Var[\epsilon_t]+1}\frac{k+1}{r - 1}
\label{eq:divided_simplified_n1}
\end{equation}
Multiply equation \eqref{eq:divided_simplified_n1} by \( k + 1 \):
\begin{equation}
- (r + 1)(k + 1)^2 + 2 r (k + 2)(k + 1) - (r - 1)(k + 2)(k + 3) = \frac{Var[\epsilon_t]}{Var[\epsilon_t]+1}\frac{(k+1)^2}{r - 1}
\label{eq:multiplied_n12}
\end{equation}
The left side of \eqref{eq:multiplied_n12} is 
$$
- r k^2 - 2 r k - r - k^2 - 2k - 1 + 2 r k^2 + 6 r k + 4 r - r k^2 - 5 r k - 6 r + k^2 + 5k + 6 
$$
Thus, with simplification of the left side of \eqref{eq:multiplied_n12}, we arrive at 
\begin{equation} 
-r(k+3)+3k+5 =  \frac{Var[\epsilon_t]}{Var[\epsilon_t]+1}\frac{(k+1)^2}{r - 1}
\label{eq:multiplied_nr}
\end{equation}
Since $r>1$, it follows that 
\begin{equation}
-r(k+3)+3k+5 \ge 0 
\label{eq:multiplied_n1}
\end{equation}
and 
\begin{equation}
1<r \leq \frac{3k+5}{k+3}<3 
\label{eq:multiplied_q}
\end{equation}
Multiply both sides of \eqref{eq:multiplied_nr} by \( (r - 1) \):
\begin{equation}
(r - 1)\Big( -r(k + 3) + 3k + 5 \Big) =  \frac{Var[\epsilon_t]}{Var[\epsilon_t]+1}(k+1)^2
\label{eq:multiplied_qq}
\end{equation}
Simplification of the left side of \eqref{eq:multiplied_qq} leads to:
\begin{equation}
-r^2(k + 3) + 4kr + 8r - 3k - 5 =  \frac{Var[\epsilon_t]}{Var[\epsilon_t]+1}(k+1)^2
\label{eq:multiplied_qq1}
\end{equation}
Now express \eqref{eq:multiplied_qq1} as a quadratic equation in terms of \( r \):
\begin{equation}
(k + 3) r^2 - \left( 4k + 8 \right) r + \left( 3k + 5 + \dfrac{Var[\epsilon_t]}{Var[\epsilon_t] + 1}(k + 1)^2 \right) = 0.
\label{eq:multiplied_qq12}
\end{equation}
Let \( A \), \( B \), and \( C \) in \eqref{eq:multiplied_qq12} for computing its discriminant as:
\begin{align*}
A = k + 3, B = -\left( 4k + 8 \right) = -4k - 8, C = 3k + 5 + \dfrac{Var[\epsilon_t]}{Var[\epsilon_t] + 1}(k + 1)^2
\end{align*}
and
\begin{align*}
B^2 &= 16k^2 + 64k + 64,\\
4AC &= 4(k + 3)\left( 3k + 5 + \dfrac{Var[\epsilon_t]}{Var[\epsilon_t] + 1}(k + 1)^2 \right) \\
&= 4(k + 3)(3k + 5) + \dfrac{4Var[\epsilon_t](k + 3)(k + 1)^2}{Var[\epsilon_t] + 1}\\
  &=12k^2 + 56k + 60 + \dfrac{4Var[\epsilon_t](k + 3)(k + 1)^2}{Var[\epsilon_t] + 1}.
\end{align*}
Thus the discriminant \( D=B^2-4AC \) of the quadratic equation with respect to $r$ is given by:
\begin{align*}
D &= \left( 16k^2 + 64k + 64 \right) - \left( 12k^2 + 56k + 60 + \dfrac{4Var[\epsilon_t](k + 3)(k + 1)^2}{Var[\epsilon_t] + 1} \right) \\
&= 4k^2 + 8k + 4 - \dfrac{4Var[\epsilon_t](k + 3)(k + 1)^2}{Var[\epsilon_t] + 1}\\
&=  4(k + 1)^2 \left( \dfrac{1 - Var[\epsilon_t](k + 2)}{Var[\epsilon_t] + 1} \right).
\end{align*}
For real solutions, the discriminant \( D \) must be non-negative. Therefore, the condition for real solutions is:
\begin{equation}
0 < k \leq \dfrac{1}{Var[\epsilon_t]}-2.
\label{eq:multiplied_q1}
\end{equation}
We now have the maximum of the feasible variance of $\epsilon_t$ 
\begin{equation}
0 \leq Var[\epsilon_t] \leq 0.5.
\label{eq:multiplied_q2}
\end{equation}
Using the quadratic formula, we have
\begin{align}
r_{\pm} &= \dfrac{ 4k + 8 \pm 2(k + 1) \sqrt{ \dfrac{1 - Var[\epsilon_t](k + 2)}{Var[\epsilon_t] + 1} } }{ 2(k + 3) } \\
&= \dfrac{ 2k + 4 \pm (k + 1) \sqrt{ \dfrac{1 - Var[\epsilon_t](k + 2)}{Var[\epsilon_t] + 1} } }{ k + 3 }.
\label{eq:multiplied_q3}
\end{align}
It is easy to see that both $r_{\pm} \geq 1$.   Since $r_{+} > r_{-}$ and $\dfrac{1 - Var[\epsilon_t](k + 2)}{Var[\epsilon_t] + 1} \leq 1$ and therefore 
$$2k + 4 - (k + 1) \sqrt{ \dfrac{1 - Var[\epsilon_t](k + 2)}{Var[\epsilon_t] + 1} } \geq  k + 3 $$
We summarize the above result as Theorem \ref{theorm1}.

\subsection{Properties of Gamma Distribution}
\subsubsection{Identity of Gamma Distribution}\label{section5_3}

In this section, our goal is to prove the following identity for $s \geq 0$ 
\begin{equation}
E\left[ X_t^n e^{- s X_t} \right] = \frac{\Gamma(k + n)}{\Gamma(k)} \frac{\theta^n}{(1 + s \theta)^{k + n}}
\label{identity}
\end{equation}
where \( X_t \sim \text{Gamma}(k, \theta) \)
\begin{equation}
f_{X_t}(x) = \frac{1}{\Gamma(k) \theta^{k}} x^{k - 1} e^{- x / \theta}, \quad x > 0.
\end{equation}

\noindent
We aim to compute the expected value:

\begin{eqnarray}
E\left[ X_t^n e^{- s X_t} \right] &=& \int_{0}^{\infty} x^{n} e^{- s x} f_{X_t}(x) \, dx\\
%
&=& \int_{0}^{\infty} x^{n} e^{- s x} \left( \frac{1}{\Gamma(k) \theta^{k}} x^{k - 1} e^{- x / \theta} 
\right) dx \\
&=& \frac{1}{\Gamma(k) \theta^{k}} \int_{0}^{\infty} x^{n + k - 1} e^{- x (s + 1 / \theta)} dx
\label{5_3_exp}
\end{eqnarray}
Recall the standard gamma integral

\begin{equation}
\int_{0}^{\infty} x^{\alpha - 1} e^{- \lambda x} dx = \frac{\Gamma(\alpha)}{\lambda^{\alpha}}, \quad \text{for } \lambda > 0, \ \alpha > 0
\label{gamma}
\end{equation}
Applying \eqref{gamma} to \eqref{5_3_exp} with $\alpha = n+k,\, \lambda = s+1/\theta$, we get \eqref{identity}
\begin{eqnarray*}
E\left[ X_t^n e^{- s X_t} \right] &=& \frac{1}{\Gamma(k) \theta^{k}} \cdot \frac{\Gamma(n + k)}{\left( s + \dfrac{1}{\theta} \right)^{n + k}} \\
%
&=& \frac{ \Gamma(k + n) }{ \Gamma(k) } \frac{ \theta^{n} }{ (1 + s \theta )^{ k + n } }
\end{eqnarray*}


\subsubsection{Third Central Moment for Gamma Distribution}
Assume that 
\[
X_n \sim \operatorname{Gamma}(k,\theta),
\]
so that
\[
E[X_n] = k\theta \quad \text{and} \quad \operatorname{Var}(X_n) = k\theta^2.
\]
Define the fluctuation
\[
\beta_n = X_n - E[X_n] = X_n - k\theta.
\]
The third central moment is given by
\[
E[\beta_n^3] = E\Bigl[(X_n - k\theta)^3\Bigr].
\]
We will derive the following estimation. 

\begin{equation}
E[\beta_n^3] = O\Bigl((\operatorname{Var}(X_n))^{3/2}\Bigr).
\label{estimation1}
\end{equation}

Express this in terms of the raw moments:
\[
E[\beta_n^3] = E[X_n^3] - 3E[X_n]\,E[X_n^2] + 2E[X_n]^3.
\]
Using the identity
\[
E[X_n^n] = \frac{\Gamma(k+n)}{\Gamma(k)}\,\theta^n,
\]
we obtain:
\[
E[X_n] = k\theta,\quad E[X_n^2] = (k+1)k\,\theta^2,\quad E[X_n^3] = (k+2)(k+1)k\,\theta^3.
\]
Substitute these into the expression for \(E[\beta_n^3]\):
\[
\begin{aligned}
E[\beta_n^3] &= (k+2)(k+1)k\,\theta^3 - 3\,(k\theta)\,( (k+1)k\,\theta^2) + 2\,(k\theta)^3 \\
&= k\,\theta^3\Bigl[(k+2)(k+1) - 3k(k+1) + 2k^2\Bigr].
\end{aligned}
\]
Notice that
\[
(k+2)(k+1) = k^2 + 3k + 2.
\]
Then,
\[
\begin{aligned}
(k+2)(k+1) - 3k(k+1) + 2k^2 
&= (k^2 + 3k + 2) - 3k^2 - 3k + 2k^2 \\
&= k^2 + 3k + 2 - 3k^2 - 3k + 2k^2 \\
&= 2.
\end{aligned}
\]
Thus,
\[
E[\beta_n^3] = k\,\theta^3\cdot 2 = 2k\,\theta^3.
\]
Since 
\[
\operatorname{Var}(X_n) = k\,\theta^2,
\]
we have
\[
(\operatorname{Var}(X_n))^{3/2} = (k\,\theta^2)^{3/2} = k^{3/2}\,\theta^3.
\]
or equivalently,
\[
E[\beta_n^3] = O\Bigl((\operatorname{Var}(X_n))^{3/2}\Bigr).
\]

\subsubsection{Fourth Central Moment for Gamma Distribution}

Assume that 
\[
X_t\sim \operatorname{Gamma}(k,\theta),
\]
so that the mean and variance are
\[
E[X_t] = k\theta,\quad \operatorname{Var}(X_t)=k\theta^2.
\]
We define the fluctuation
\[
\beta_t = X_t - E[X_t].
\]

Our goal is to show that

\begin{equation}
E[\beta_t^4] = O\Bigl((\operatorname{Var}(X_t))^{2}\Bigr).
\label{estimation2}
\end{equation}
The identity \eqref{identity}
implies that by setting \(s=0\) we obtain the raw moments:
\[
E[X_t^n] = \frac{\Gamma(k+n)}{\Gamma(k)}\,\theta^n.
\]
In particular, we have:
\begin{align*}
E[X_t] &= \frac{\Gamma(k+1)}{\Gamma(k)}\,\theta = k\,\theta,\\[1mm]
E[X_t^2] &= \frac{\Gamma(k+2)}{\Gamma(k)}\,\theta^2 = (k+1)k\,\theta^2,\\[1mm]
E[X_t^3] &= \frac{\Gamma(k+3)}{\Gamma(k)}\,\theta^3 = (k+2)(k+1)k\,\theta^3,\\[1mm]
E[X_t^4] &= \frac{\Gamma(k+4)}{\Gamma(k)}\,\theta^4 = (k+3)(k+2)(k+1)k\,\theta^4.
\end{align*}

The fourth central moment (which is \(E[\beta_t^4]\)) is given by:
\[
E[(X_t-E[X_t])^4] = E[X_t^4] - 4E[X_t]\,E[X_t^3] + 6\bigl(E[X_t]\bigr)^2E[X_t^2] - 3\bigl(E[X_t]\bigr)^4.
\]
Substitute the moments computed above:
\begin{align*}
E[\beta_t^4] &= (k+3)(k+2)(k+1)k\,\theta^4 \\
&\quad - 4\,(k\,\theta)\,(k+2)(k+1)k\,\theta^3 \\
&\quad + 6\,(k\,\theta)^2\,(k+1)k\,\theta^2 \\
&\quad - 3\,(k\,\theta)^4.
\end{align*}
Factor out the common factor \(k^2\theta^4\) and simplify them:  
\[
E[\beta_t^4] = \left(3+\frac{6}{k}\right)k^2\theta^4,
\]
Since
\[
\operatorname{Var}(X_t) = k\theta^2,
\]
we have
\[
\operatorname{Var}(X_t)^2 = k^2\theta^4.
\]
Therefore, we conclude the verification of \eqref{estimation2}.

\section{Acknowledge}
This work by YW was partially supported by NSF-DMS grant 1752709.


\begin{thebibliography}{9}


\bibitem{may1976simple} R. M. May, \textit{Simple mathematical models with very complicated dynamics}, Nature, vol. 261, no. 5560, pp. 459--467, 1976.

\bibitem{feigenbaum1980quantitative} M. J. Feigenbaum, \textit{Quantitative universality for a class of nonlinear transformations}, Journal of Statistical Physics, vol. 19, no. 1, pp. 25--52, 1978.

\bibitem{hilborn2000chaos} R. C. Hilborn, \textit{Chaos and Nonlinear Dynamics: An Introduction for Scientists and Engineers}, 2nd ed., Oxford University Press, 2000.


%
%

\bibitem{hurst1951long} H. E. Hurst, \textit{Long-term storage capacity of reservoirs}, Transactions of the American Society of Civil Engineers, vol. 116, pp. 770--799, 1951.

\bibitem{eckmann1987recurrence} J.-P. Eckmann, S. O. Kamphorst, and D. Ruelle, \textit{Recurrence plots of dynamical systems}, Europhysics Letters, vol. 4, no. 9, pp. 973--977, 1987.

\bibitem{Allen}
Allen, L. J. S. (2010). \emph{An Introduction to Stochastic Processes with Applications to Biology}. Chapman and Hall/CRC.

\bibitem{Renshaw}
Renshaw, E. (1991). \emph{Modelling Biological Populations in Space and Time}. Cambridge University Press.

\bibitem{May2001}
R. M. May, (1974) \emph{Stability and Complexity in Model Ecosystems}, Princeton U.P., Princeton.

\bibitem{Dennis1984}
B. Dennis and G. P. Patil, \emph{The Gamma Distribution and Weighted Multimodal Gamma Distributions as Models of Population Abundance}, Mathematical Biosciences, 68, 187–212, 1984.

\bibitem{Peters1989}
C. S. Peters, M. Mangel, R. F. Costantino, \emph{Stationary Distribution of Population Size in Tribolium}, Bulletin of Mathematical Biology, 51(5), 625–638, 1989.

\bibitem{DennisCostantino1988}
B. Dennis and R. F. Costantino, \emph{Analysis of Steady-State Populations with the Gamma Abundance Model and its Application to Tribolium}, Ecology, 69(4), 1200–1213, 1988.

\bibitem{Costantino1981}
R. F. Costantino and R. A. Desharnais, \emph{Gamma Distributions of Adult Numbers for Tribolium Populations in the Region of their Steady-States}, Journal of Animal Ecology, 50, 667–681, 1981.


\bibitem{Matis2003}
James H. Matis, Thomas R. Kiffe b, Eric Renshawc, Janet Hassan, \emph{A simple saddlepoint approximation for the equilibrium
distribution of the stochastic logistic model of population growth}, Ecological Modelling, 161 (2003) 239–248.


\bibitem{wikipedia_gamma}
Wikipedia contributors, "Gamma distribution," \emph{Wikipedia, The Free Encyclopedia}, 2024.


\bibitem{kot2001elements}
M. Kot, \emph{Elements of Mathematical Ecology}. Cambridge University Press, 2001.

%

\bibitem{may1976simple} May, R. M. (1976). \emph{Simple mathematical models with very complicated dynamics.} \emph{Nature}, 261(5560), 459-467.

\bibitem{strogatz2018nonlinear} Strogatz, S. H. (2018). \emph{Nonlinear Dynamics and Chaos: With Applications to Physics, Biology, Chemistry, and Engineering}. CRC Press.




\bibitem{Yan2024} Dingding Yan, Mengqi He, Robert A. Cheke, Qianqian Zhang, Sanyi Tang, (2024). \emph{A stochastic hormesis Ricker model and its application to multiple fields}. Chaos, Solitons \& Fractals, Volume 185, 2024, 115131

\bibitem{Nasell2003} Ingemar Nåsell,(2003). \emph{Moment closure and the stochastic logistic model}. Theoretical Population Biology, Volume 63, Issue 2, 2003, Pages 159-168,

\bibitem{Pielou1975} E. C. Pielou (1975). \emph{Ecological Dirersity},  Wiley, New York, 1975.

\bibitem{Engen1978} S. Engen, (1978). \emph{ Stochustic Abundunce Models}, Chapman and Hall, London, 1978.

\bibitem{Schreiber2021} Schreiber, Sebastian J. and Huang, Shuo and Jiang, Jifa and Wang, Hao (2021). \emph{Extinction and Quasi-Stationarity for Discrete-Time, Endemic SIS and SIR Models}, SIAM Journal on Applied Mathematics, 81(5) 2195-2217.


\bibitem{earlywarning2023} Lingyu Feng, Ting Gao, Wang Xiao, Jinqiao Duan,  \textit{Early Warning Indicators via Latent Stochastic Dynamical Systems}, Chaos: An Interdisciplinary Journal of Nonlinear Science, vol. 34, pp. 031101, 2024. 
    
\bibitem{arxiv2019chaos} Daniel Toker, Friedrich T. Sommer,  Mark D’Esposito, \textit{A Simple Method for Detecting Chaos in Nature}, Communications Biology volume 3, Article number: 11 (2020)
     
\bibitem{chaos2017extremes} C. Franzke  \textit{Extremes in Dynamic-Stochastic Systems}, Chaos: An Interdisciplinary Journal of Nonlinear Science,   vol. 27, no. 1, pp. 012101, 2017. 
%
\bibitem{physrev2022lineagetrees} Hagai Rappeport, Irit Levin Reisman, Naftali Tishby, Nathalie Q. Balaban \textit{Detecting Chaos in Lineage-Trees: A Deep Learning Approach}, vol. 4, no. 1, pp. 013223, 2022.
%
%
%
%
%
\bibitem{wolf1985lyapunov} A. Wolf, J. B. Swift, H. L. Swinney, and J. A. Vastano, \textit{Determining Lyapunov exponents from a time series}, Physica D: Nonlinear Phenomena, vol. 16, no. 3, pp. 285--317, 1985. 
%
\bibitem{may1976simple} R. M. May, \textit{Simple mathematical models with very complicated dynamics}, Nature, vol. 261, no. 5560, pp. 459--467, 1976.
%
%
\bibitem{ricker1954stock} W. E. Ricker, \textit{Stock and recruitment}, Journal of the Fisheries Research Board of Canada, vol. 11, no. 5, pp. 559--623, 1954.
%
\bibitem{may1974biological} R. M. May, \textit{Biological populations with nonoverlapping generations: stable points, stable cycles, and chaos}, Science, vol. 186, no. 4164, pp. 645--647, 1974.
%
\bibitem{moss1989noise} F. Moss and P. V. E. McClintock, \textit{Noise in Nonlinear Dynamical Systems}, vol. 3, Cambridge University Press, 1989.
%
%
%
%
%
%
%
%
%
%
%
\bibitem{crutchfield1982fluctuations} 
J. P. Crutchfield, J. D. Farmer, and B. A. Huberman, 
\textit{Fluctuations and Simple Chaotic Dynamics}, 
Physics Reports, vol. 92, no. 2, pp. 45--82, 1982.

%
\bibitem{erguler2008statistical} 
K. Erguler and M. P. H. Stumpf, 
\textit{Statistical Interpretation of the Interplay Between Noise and Chaos in the Stochastic Logistic Map}, 
Mathematical Biosciences, vol. 216, pp. 90--99, 2008.
%
\bibitem{hamzi2024simplicity} 
B. Hamzi and K. Dingle, 
\textit{Simplicity Bias, Algorithmic Probability, and the Random Logistic Map}, 
Physica D: Nonlinear Phenomena, vol. 463, p. 134160, 2024.
%
\bibitem{li2008effects} 
F.-g. Li, 
\textit{Effects of Noise on Periodic Orbits of the Logistic Map}, 
Central European Journal of Physics, vol. 6, no. 3, pp. 539--545, 2008.
%
\bibitem{sato2018dynamical} 
Y. Sato, T. S. Doan, J. S. W. Lamb, and M. Rasmussen, 
\textit{Dynamical Characterization of Stochastic Bifurcations in a Random Logistic Map}, 
arXiv preprint, arXiv:1811.03994, 2018.

\bibitem{li1975} 
Tien-Yien Li; James A. Yorke, 
\textit{Period Three Implies Chaos}, 
The American Mathematical Monthly, Vol. 82, No. 10. (Dec., 1975), pp. 985-992.



%
%
%
%
%
%
%
%
%
%
%
%
%
%
%

\bibitem{wang2025}
Haiyan Wang,
\textit{Equilibrium analysis of discrete stochastic population models with gamma distribution},
Mathematical Biosciences,  381, 109398, 2025, https://doi.org/10.1016/j.mbs.2025.109398

\end{thebibliography}
\end{document}